\input amstex
\input amsppt1
\loadmsbm
\UseAMSsymbols
\NoBlackBoxes
\parindent=0pt
\magnification=1200

\font\pio=cmr10 scaled 1230
\font\pioa=cmb10 scaled 1300
\def\empty{{\text{\rm \pio \o}}}

\def\cd{\Cal D}
\def\td{\widetilde{\Cal D}}

\def\ca{{\Cal A}}

\def\cs{\Cal S}

\def\enpr{\quad \vrule height .9ex width .8ex depth -.1ex}

\def\bar{\overline}

\def\bre{\Bbb R}
\def\bna{\Bbb N}

\def\bte{\bar\theta}

\def\llim{\mathop{\longrightarrow}}
\def\vo{V^{(0)}}
\def\von{V^{(1)}}
\def\rvo{\bre^{V^{(0)}}}
\def\rvon{\bre^{V^{(1)}}}

\def\Lb{\Lambda}
\def\cl{{\Cal L}}
\def\ce{{\Cal E}}

\def\cf{{\Cal F}}
\def\cg{{\Cal G}}

\def\qua{\quad $\,$}
\def\smad{\smallskip Proof.\ }
\def\witi{\widetilde}
\def\disp{\displaystyle}
\def\s,{\quad $\,$}

\font\pio=cmr10 scaled 1230
\def\empty{{\text{\rm \pio \o}}}

\def\enpr{\quad \vrule height .9ex width .8ex depth -.1ex}

\def\bar{\overline}

\def\bre{\Bbb R}
\def\bna{\Bbb N}

\def\defin#1.#2{\noindent \bigskip {${\bf Definition \ #1.#2}$}   }
\def\osser#1.#2{\noindent \bigskip {${\bf Remark     \ #1.#2}$}   }
\def\propo#1.#2{\noindent \bigskip {${\bf Proposition\ #1.#2}$}   }
\def\teore#1.#2{\noindent \bigskip {${\bf Theorem    \ #1.#2}$}   }
 \def\corol#1.#2{\noindent \bigskip {${\bf Corollary  \ #1.#2}$}   }
\def\lemma#1.#2{\noindent \bigskip {${\bf Lemma      \ #1.#2}$}   }

\def\veps{\varepsilon}
\def\osc{ \text{\rm Osc}}

\def\llim{\mathop{\longrightarrow}}
\def\vo{V^{(0)}}
\def\von{V^{(1)}}
\def\rvo{\bre^{V^{(0)}}}
\def\rvon{\bre^{V^{(1)}}}
\def\vin#1{V^{(#1)}}
\def\rvin#1{\bre^{V^{(#1)}}}

\def\Lb{\Lambda}
\def\cl{{\Cal L}}
\def\ce{{\Cal E}}
\def\cf{{\Cal F}}
\def\rcf{{\bre^{\Cal F}}}
\def\cg{{\Cal G}}

\def\qua{\quad $\,$}
\def\smad{\smallskip Proof.\ }
\def\witi{\widetilde}
\def\disp{\displaystyle}
\def\s,{\quad $\,$}

\font\pio=cmr10 scaled 1230
\def\empty{{\text{\rm \pio \o}}}
 
\def\cite#1{[#1]}

 \centerline  {\bf \pioa Integrals of convex functions in the gradients on fractals.}
\medskip\centerline{\bf Roberto Peirone}

\smallskip\centerline{\it Universit\`a di Roma-Tor Vergata,Dipartimento di Matematica}          
\centerline{\it Via della Ricerca Scientifica, 00133 Roma, Italy}
\smallskip
\centerline{E.mail: peirone\@ mat.uniroma2.it}

\bigskip\medskip

{\bf Abstract}. {\it In this paper, I describe the construction of certain functional
integrals in the gradient on finitely ramified fractals, which have 
a sort of self-similarity property.}

\bigskip\medskip

\centerline{\bf 1. Introduction.}

\medskip 
\qua 
The subject of this paper is analysis on finitely ramified 
self-similar fractals. 
 Some examples of self-similar fractals are
  the (Sierpinski) Gasket, the Vicsek set and the (Lindstr\o m) Snowflake.
  The Gasket is obtained starting with an equilateral triangle, next
  dividing it into four triangles  having edge equal to one half of the original triangle
  and removing the central triangle, next repeating the same construction 
  on every of the remaining triangles and so on.
  The Vicsek set  is constructed starting with a square and dividing it into nine squares
  whose edges are ${1\over 3}$ of the original edge, and 
 considering  only
  the four squares at the vertices and the central square. 
  The  Snowflake instead is constructed 
  starting with a regular
   hexagon and dividing it into seven hexagons, the relative contractions 
   having for fixed points the
   six vertices and the centre of the hexagon. 
   A general way to construct  a
 (self-similar) fractal is the following.  We start with finitely many contractive 
 (i.e., having factor $<1$) similarities 
 $\psi_1,...,\psi_k$ in $\bre^{\nu}$ (more generally, in a compact metric space). 
 Then, the self-similar fractal generated by such similarities is 
 the (unique) nonempty compact set $\cf$ in ${\bre}^{\nu}$ such that

$$\cf=\bigcup\limits_{i=1}^k \psi_i(\cf), \eqno (1.1)$$
 For example, in the Gasket the maps
   $\psi_i$, $i=1,2,3$, are the rotation-free contractions 
 with factor ${1\over 2}$ that have as fixed points $P_i$, the three vertices of the
 triangle,  in formula 
  $\psi_i(x)=P_i+{1\over 2}(x-P_i)$.   We also require that the fractals 
   are connected (thus excluding the Cantor set) and finitely ramified,
   which more or less means that the copies of the fractal via the similarities
   intersect only at finitely many points. An example of fractal
  that is not finitely ramified is the Sierpinski Carpet.
  Note also that the segment-line $[0,1]$ can be seen
  as a (degenerate) finitely ramified self-similar fractal with two maps
  $\psi_1$ and $\psi_2$ defined by $\psi_1(x)={x\over 2}$,
  $\psi_2(x)={1+x\over 2}$.

  \qua In the present paper, I will consider  the
 P.C.F. self-similar sets, a class of 
 finitely ramified fractals introduced by J. Kigami  in \cite{3}, with a mild additional requirement
 as in \cite{1} or \cite{6}.
  A description
    of the general theory of P.C.F. self-similar sets with many examples 
    (including those described above) can be 
    found in \cite{4}.  
    
  \qua  One of the
  problems widely discussed in this area is  the construction 
  of self-similar irreducible
  Dirichlet forms, that is analogs of the Dirichlet integral, on
 fractals. This can be also interpreted as the construction 
  of diffusions on fractals. 
   A Dirichlet form  $\ce$ on the fractal $\cf$
 is a functional from $\bre^{\cf}$ to $[0,+\infty]$
which is a quadratic form and satisfies
$\ce(v+c)=\ce(v)$ for every constant $c$ and certain additional
properties. We say that it is irreducible if  it takes the value $0$ only 
at the constants, and we say that it is
 self-similar  if  satisfies
 
 $${1\over\rho}\sum\limits_{i=1}^k \ce(v\circ\psi_i)= \ce(v) \eqno(1.2)$$
 
for some positive $\rho$. 
%In particular, when $K$ is the segment-line, we have $\sigma={1\over 2}$
%and $\theta=2$.
In general finitely ramified fractals the existence of a self-similar energy is
a delicate problem which reduces to an eigenvector problem
 of a special nonlinear operator $\Lb$ defined on the set of the  Dirichlet forms 
 on a specific {\it finite subset}  $\vo$ of the
 fractal, which can be interpreted as the boundary fractal.
 In fact, $\vo$ is  a subset of the set of the fixed points of $\psi_i$, $i=1,...,k$, when 
 $\cf$ is defined via (1.1), usually formed by those fixed points that are extremum points of the convex hull
 of the set of all fixed points. The energy $\ce$ on $\cf$ is defined based on the eigenvector
 $E$ of $\Lb$, which is quadratic form on $\vo$. More precisely,
 given $E$ as above, then we 
 can construct $\ce$ irreducible and self-similar such that the following holds.

 \qua {\it We have $E(u)=\min\ce(v)$ where the minimum is taken over 
the functions defined from $\cf$ with values into $\bre$
 that amount to $u$ on $\vo$. Moreover the function $v$ realizing the minimum attains
 its maximum and its minimum on $\vo$.}
 
Conditions for the existence of self-similar energies are known, and in particular
this occurs for nested fractals, a class of highly symmetric fractals introduced by Tom 
Linsdtr\o m in \cite{5}. However, in a lot of fractals such a self-similar energy does not exist.
It appeared to be natural to consider a variant of self-similarity formula
(1.1), placing {\it weights} on the cells.
The existence of a self-similar energy in this broader sense 
on P.C.F. self-similar sets is still an open problem.
However, in \cite{6} it is proved that this occurs for fractals with connected interior,
and that, in general fractals, this occurs at a suitable level (depending on the fractal).
We can proceed similarly when $\ce$ is not  a quadratic form but is $p$-homogeneous.
In \cite{2} the existence of a $p$-homogeneous energy
is proved on the Gasket and, more generally, on fractals having special
symmetry property (usually stronger than the symmetry
property of nested fractals) called {\it weakly completely symmetric}.

\qua In the present paper, I investigate the existence of self-similar 
energies  when more generally, $\ce$
 is a convex functional having certain additional properties (see Section 3).  
 In other words, I discuss functionals on fractals, analogous to 
 
 $$\ce(v)=\int\limits_A F(\nabla v)\eqno (1.3)$$
 
  in regions $A$ contained in euclidean spaces
 when $F$ is a convex function.
 Note that, thanks to the $2$-homogeneity,
 (1.2) can be interpreted
alternatively as

$$\sum\limits_{i=1}^k \ce(\theta v\circ\psi_i)= \ce(v) \eqno(1.4)$$

or also

$${1\over \sigma}\sum\limits_{i=1}^k \ce(\theta v\circ\psi_i)= \ce(v) \eqno(1.5)$$

for some positive $\sigma$ and  $\theta$, satisfying
 ${\theta^2\over\sigma}={1\over\rho}$.
Formulas (1.4) and (1.5), are, in my view, more natural, in that they are more related to
the notion of a functional in the derivative. Moreover, they appear to be more
appropriate for the more general case of a convex functional. 
In fact, in the case $\cf$ is the segment-line, then
for the functional described in (1.3) with $A=[0,1]$,
 (1.5) holds with $\sigma=2$ and $\theta=2$.
In the case of general fractals it appears to be too restrictive to require that
$\sigma $ and $\theta$ are constant. Thus, I will prove that (1.5) holds 
for fixed $\sigma\in]0,1]$
and $\theta$ continuously depending on the restriction of $v$ to $\vo$. Note that 
by continuity, (1.5) is easily seen to
 hold for some $\theta$ depending on $v$, and it is instead a
 nontrivial  requirement that
$\theta$ only depends on the restriction of $v$ to $\vo$.
I will prove, in fact, a 
 more precise statement (see Theorems 6.1 and 6.2) which is in some sense, an analog 
 to the statement for quadratic form. Namely,
suppose given a convex functional $E$ (not necessarily a quadratic form) from
$\rvo$ to $[0,+\infty[$ having certain additional properties, which are described 
in Section 3. One of them, in particular, states that
$E$ approximates an eigenform when the function $u$ on $\rvo$ tends to $0$.
Then, we  can associate 
to $E$ a unique semicontinuous functional $\ce$ from $\bre^{\cf}$
to $\bre$ such that

\qua a) $\ce$ is self-similar in the sense of (1.5) with $\theta$
continuously depending on the restriction of $v$ to $\vo$ 

\qua b) $E(u)=\min\ce(v)$ where the minimum is taken over 
the functions defined from $\cf$ with values into $\bre$
 that amount to $u$ on $\vo$. Moreover the function $v$ realizing the minimum attains
 its maximum and its minimum on $\vo$.

\qua In the case of the segment-line $[0,1]$ the functional $\ce$ defined in
(1.3) is associated to $E$ defined by $E(u)=F\big(v(1)-v(0)\big)$.

 \bigskip

\bigskip

\centerline{\bf 2. Notation.}
\bigskip

In the present Section,  I fix the general setting and
give the preliminary results.
 In view of  (1.1),   we can define a fractal  by giving 
 a finite set $\Psi$ of contractions 
 $\psi_i$. In order to define the fractal,
 I here follow an approach  similar
 to that introduced    in   \cite{1}  and  already
 discussed in previous papers of mine (see e.g., \cite{6}).
Let $\Psi=\big\{\psi_1,...,\psi_k\big\}$  be a set
     of one-to-one 
   maps defined on a finite set $V=\vo=\big\{P_1,...,P_N\big\}$ 
   (not necessarily a subset of 
   $\bre^{\nu}$),  with $2\le N\le k$, and put
   
   $$\von=\bigcup\limits_{\psi\in\Psi}\psi(V^{(0)})\, .\eqno (2.1)$$ 
   
   We  require that for each $j=1,...,N$  
   
   $$\psi_j(P_j)=P_j,\qquad  P_j\notin\psi_i(V^{(0)})\quad\forall\,i\ne j
   \ \forall\,j=1,...,N\eqno (2.2)$$
   $$  \von \text{\rm\ is\ connected.}\eqno (2.3)$$

  Here, we say that $\von$ is connected if  for every $i,i'=1,...,k$
  we can find a sequence of indices $i_0,...,i_h=1,...,k$ such that
  $i_0=i$, $i_h=i'$, and
  $V_{i_{s-1}}\cap V_{i_s}\ne\empty$ for every $s=1,...,h$,
   where   $V_i:=\psi_i(\vo)$, $i=1,...,k$.

    \s, Given a set $\Psi$ of maps satisfying the previous properties, 
    we can  construct a self-similar finitely ramified fractal $\cf(\Psi)$ or simply $\cf$ 
    (embedded in a metric space) having $\Psi$ (more precisely a set of maps 
    whose restrictions to $\vo$ form $\Psi$) as the set of contractions.
    We can in fact construct 
    different fractals with that property, but they are "isomorphic". 
    More precisely, $\cf$ has the following properties.
    \smallskip
    
    \qua $P_1$) $\cf$ is a compact and connected metric space with distance $d$ containing $\vo$

     \qua $P_2$) We can extend $\psi_i$ as one-to-one maps from $\cf$ into itself  
    such that $\cf=\bigcup\limits_{i=1}^k \psi_i(\cf)$. 
   
   \qua $P_3$) Every $\psi_i$ is a similarity on $\cf$ which we can suppose to have as factor 
   ${1\over 2}$. More precisely
   $$d\big(\psi_i(x),\psi_i(y)\big)={1\over 2} d(x,y)\quad \forall x,y\in \cf\ \ \forall i=1,...,k.$$

   \qua $P_4$) We have $\psi_{i_1,...,i_n}(\cf)\cap\psi_{i'_1,...,i'_n}(\cf)=
   \psi_{i_1,...,i_n}(\vo)\cap\psi_{i'_1,...,i'_n}(\vo)$ if $(i_1,...,i_n)\ne (i'_1,...,i'_n)$.
   
   \smallskip
   
   \qua 
   Here,  we define $\psi_{i_1,...,i_n}$  on $\cf$ by $\psi_{i_1,...,i_n}:=
   \psi_{i_1}\circ\cdots\circ\psi_{i_n}$ and we call {\it $n$-cells} the sets
   $V_{i_1,...,i_n}=\psi_{i_1,...,i_n}(\vo)$
   and
   more generally, we define 
   $$A_{i_1,...,i_n}:=\psi_{i_1,...,i_n}(A)$$ 
   for every subset $A$ of $\cf$.
   Moreover, we put $\vin n=\bigcup\limits_{i_1,...,i_n=1}^k V_{i_1,...,i_n}$.
   It can be easily proved that the sequence of sets $\vin n$ is increasing in $n$, and that
   its union $\vin{\infty}$ is dense in $\cf$.
  We can characterize the continuity of a function from $\cf$ to $\bre$ or the
  uniform convergence of $v_n$ to $v$ on $\cf$ in terms of the behavior on  
  the sets $\cf_{i_1,...,i_n}$. Given a function
  $f$ from a set $A$ into $\bre$, we define $\osc(f)=\sup\{f(x)-f(y):x,y\in A\}$
and $\osc_B (f)=\osc(f|_B)$ whenever $B\subseteq A$. 
  Then

  \bigskip
  {\bf Lemma 2.1.} {\sl We have
  
\qua   i) A function $v\in\rvin {\infty}$ is uniformly continuous if and only if
$$\max\big\{\osc(v\circ\psi_{i_1,...,i_n}): i_1,...,i_n=1,...,k\big\}
\llim\limits_{n\to +\infty} 0.$$
  
 \qua  ii) If $v_n,v\in \rcf$, then $v_n\llim\limits_{n\to +\infty} v$ 
 uniformly  if and only if
$$\max\limits_{i_1,...,i_n=1,...,k}\big\{\sup\limits_{\cf_{i_1,...,i_n}}:|v_n-v|\big\}
\llim\limits_{n\to +\infty} 0.$$}
\smad
The proof of i) is standard (see for example [5]). The idea is that 

$\max\limits_{i_1,...,i_n=1,...,k}$diam$(\cf_{i_1,...i_n})\llim\limits_{n\to +\infty}0$, combined with 
the fact that, by $P_4)$, 

$${\vin \infty}_{i_1,...,i_n}\cap {\vin \infty}_{i'_1,...,i'_n}=\empty \iff
{\cf}_{i_1,...,i_n}\cap {\cf}_{i'_1,...,i'_n}=\empty.$$
The proof of ii) is trivial as $\cf=\bigcup\limits_{i_1,...,i_n=1,...,k}
\cf_{i_1,...,i_n}$.
\enpr

 \bigskip 
  \qua I now recall the notion  of  Dirichlet forms
  on $\vo$, and the renormalization operator defined on it. 
  Let  $J=\{\{j_1,j_2\}:j_1,j_2=1,...,N, j_1\ne j_2\}$.
  I will denote by
$\cd(\vo)$ or simply $\cd$ the set
of the  Dirichlet forms on $\vo$, invariant with respect to
an additive constant,  i.e., the set of the functionals $E$ 
from $\rvo$ into $\bre$ of the form

 $$E(u)=\sum\limits_{\{j_1,j_2\}\in J}
  c_{\{j_1,j_2\}}(E)\big(u(P_{j_1})-u(P_{j_2})\big)^2 $$
 
where the {\it coefficients} $c_{\{j_1,j_2\}}(E)$ (or simply 
$c_{\{j_1,j_2\}}$)   of $E$ are required to
be nonnegative. I will denote by
$\td(V)$ or simply $\td$ the set of the irreducible Dirichlet forms, i.e., $E\in\td$ if $E\in
\cd$ and moreover  $E(u)=0$
 if and only if $u$ is constant.

\qua   Now,   the renormalization operator $\Lb$ is defined as follows.
  For every $u\in\rvo$  and every $E\in\td$, put
 
 $$\Lb(E)(u)=\inf\Big\{S_{1}(E)(v), v\in\cl(u)\Big\},$$
 $$S_{1}(E)(v):=\sum\limits_{i=1}^k  E(v\circ\psi_{i}),\quad 
\cl(u):= \big\{v\in \bre^{V^{(1)}}: v=u
  \ \text{on}\ V^{(0)}\big\}.$$
  
  An eigenform is an element $E$ of $\td$
   such that $\Lb(E)=\rho E$ for some positive $\rho$, which is called
  {\it eigenvalue of $E$}.
  Given an eigenform $E$ with eigenvalue $\rho$, we can 
  associate an "energy" $\ce$ on $\cf$ in the following way.
  For every $n\in\bna$ let $S_n(E)$ be defined as
  
  $$S_n(E)(v)=\sum\limits_{i_1,...,i_n=1}^k E(v\circ \psi_{i_1,...,i_n}),
  \quad \witi S_n(E)={1\over\rho^n}S_n(E).$$
  
for every $v\in\rvin n$. So, if $v\in\bre^{\cf}$, it can be easily proved that
the sequence $\witi S_n(E)(v)$ is increasing, thus tends to a
 (possibly infinite) limit which I denote by $\ce(v)$.

\bigskip
\centerline{\bf 3. Energies on $\vo$.}
\medskip

In this Section, I will generalize the previous definitions from the case of
Dirichlet forms to that of certain classes of convex functionals.
First of all, we need some generalizations of the definition 
of $\cl(u)$.
For every $u\in\rvo$, $\theta>0$, and $E:\rvo\to[0,+\infty[$ let

$$\cl'(u)= \big\{v\in \cl(u): \min u\le v\le \max u\big\}$$

Note that if $v\in\cl'(u)$  we have $\osc(v\circ\psi_i)\le \osc(u)$  for every $i=1,...,k$. So, let

$$\cl''(u)=\big\{v\in\cl'(u):\osc(v\circ\psi_i)< \osc(u)\ \ \forall\,=1,...,k\big\},$$
$$\cl'''(u)=\cl'(u)\setminus\cl''(u).$$
$$\cl(\cf,u)=\big\{v\in\bre^{\cf}: v=u \ \text{on}\ \vo\big\},$$
$$\cl'(\cf, u)= \big\{v\in \cl(\cf,u): \min u\le v\le \max u\big\}$$

where $\bre^{\cf}$ denotes the set of the {\it continuous} functions
from $\cf$ into $\bre$.
We now define analogs of $S_1(E)$ and $\Lb(E)$, but depending on $\theta$. Let

$$S_{(\theta)}(E)(v)=\sum\limits_{i=1}^k E(\theta v\circ\psi_i)\quad \forall\, v\in\rvon,$$
$$\Lb_{(\theta)}(E)(u)=\inf\big\{ S_{(\theta)}(E)(v): v\in\cl(u)\big\}.$$

We will see that, when $E$ satisfies reasonable conditions,
 the infimum in the definition of $\Lb_{(\theta)}(E)(u)$ is in fact a minimum,
but unlike the case of Dirichglet forms, there is no reason in general that it is taken at a unique point.
We denote by $H_{E,(\theta)}(u)$ the set of $v\in\cl(u)$ 
such that $\Lb_{(\theta)}(E)(u)=S_{(\theta)}(E)(v)$ and put
$$H'_{E,(\theta)}(u)=H_{E,(\theta)}(u)\cap\cl'(u).$$

I will now introduce some classes of functionals $E$ suitable for our aims. I first require some minimal properties
which guarantee some general results.
Let  $\ca_1$ be the set of $E:\rvo\to [0,+\infty[$  such that for every $u\in\rvo$ we have

\quad $Q_1$) $E$ is convex

\quad $Q_2$) $E(\pm u+c)=E(u)$ for every constant $c$.

\quad $Q_3$) $E(u)=0 \iff u$ is constant.

\quad $Q_4$) $E\big( (u\wedge a)\vee b\big)\le E(u)$ if $a\ge b$ 
and the inequality is strict if $ (u\wedge a)\vee b\ne u$.

Note that $Q_4$) is an analog of the well-known Markov property and in $\cd$ is related to the requirement
that $c_{\{j_1,j_2\}}\ge 0$. I now state some standard properties of $S_{(\theta)}(E)$
and of $\Lb_{(\theta)}(E)$.

\bigskip
{\bf Lemma 3.1.} {\sl We have

\qua i) The map $(\theta,v)\mapsto S_{(\theta)}(E)(v)$ is continuous from $]0,+\infty[\times
\rvon$ to $\bre$.

\qua ii) The map $(\theta,u)\mapsto \Lb_{(\theta)}(E)(u)$ is continuous from $]0,+\infty[\times
\rvo$ to $\bre$.

\qua iii) If $E\in\ca_1$, then for every $u\in\rvo$ and every constant $c$ we have
$\Lb_{(\theta)}(u)=\Lb_{(\theta)}(u+c)$.}
\smad
The proofs of i) and iii) are almost trivial, and the proof of ii) is  a standard argument
 (see for example Lemma 2.4 ii) in [ ] for a proof in a similar case). \enpr

\bigskip
{\bf Lemma 3.2.} {\sl For every $E\in\ca_1$, every $u\in\rvo$  and every $t\ge 1$, we have
$E(tu)\ge t E(u)$.}
\smad
This immediately follows from the convexity of $E$ and the condition $E(0)=0$. \enpr

\bigskip
{\bf Lemma 3.3.} {\sl For every $E\in\ca_1$ there exists $c>0$ such that 
$E(u)\ge c\, \osc(u)$ if $u\in\rvo$ and $\osc(u)\ge 1$.}
\smad 
Let $u$ be a nonconstant element of $\rvo$. Then,
 there exist $j_1$ and $j_2$ such that $\osc(u)=|u(P_{j_2})-u(P_{j_1})|$. Hence,
$||u-u(P_{j_1})||\ge\osc(u)$ and, also using Lemma 3.2 we have

$$E(u)=E\big(u-u(P_{j_1})\big)\ge\osc(u) E
\Big( {u-u(P_{j_1})\over || u-u(P_{j_1}) ||}\Big)\ge c\, \osc(u)\ \ \ \text{if}\ \osc(u)\ge 1$$

where $c:=\min\limits_{\cs_{j_1}}E$, and $\cs_{j_1}:= \big\{ u\in\rvo:u(P_1)=0, ||u||=1\big\}$. \enpr

\bigskip
{\bf Lemma 3.4.} {\sl For every $E\in\ca_1$, $\theta>0$ and $u\in\rvo$ the sets $H_{E,(\theta)}(u)$ 
and $H'_{E,(\theta)}(u)$ are nonempty compact and convex.}
\smad 
As $\von$ is connected, we can find  a positive integer $\bar n$ 
having the following property: 

\qua {\it for every $Q,Q'\in\von$  there exists $\bar i=1,...,k$
 such that   

$$\disp{\osc_{V_{\bar i}}(v)\ge {|v(Q)-v(Q')|\over \bar n}}.\eqno (3.1)$$}

 Thus, in view of Lemma 3.3 we have

$$S_{(\theta)}(E)(v)\llim\limits_{v\to\infty}+\infty\  \text{on} \ \cl(u),\eqno (3.2)$$

 so that,
$S_{(\theta)}(E)$ being continuous, it has a minimum on the closed set $\cl(u)$, in other words
$H_{E,(\theta)}(u)$ is nonempty. 
The compactness  of $H_{E,(\theta)}(u)$ 
easily follows from (3.2), and the convexity of $H_{E,(\theta)}(u)$
easily follows from the convexity of $S_{(\theta)}(E)$ which in turn
 immediately follows from the convexity of $E$.
We deduce from this  also the compactness and the 
convexity of $H'_{E,(\theta)}(u)$ as also $\cl'(u)$ is compact and convex. 
Finally, by $Q_4$), if $w\in H_{E,(\theta)}(u)$, then
$(w\wedge \max u)\vee \min u\in H'_{E,(\theta)}(u)$, so that $H'_{E,(\theta)}(u)$
 is nonempty. \enpr

\bigskip
{\bf Lemma 3.5.} {\sl 
For every $E\in\ca_1$, every $\sigma>0$ and every nonconstant $u\in\rvo$, there 
exists a unique $\bar\theta_{\sigma,E}(u)>0$ such that 
$\Lb_{(\bar\theta_{\sigma,E}(u))}(E)(u)=\sigma E(u)$.
Moreover, for every $\sigma>0$ and $E\in\ca$,  the map 
$\bar\theta_{\sigma,E}$ is continuous from the set of nonconstant
$u\in\rvo$ into $]0,+\infty[$.}
\smad
By Lemma 3.3 we have $E(\theta u)\le \theta E(u)$ for every $u\in\rvo$
and every $\theta\in]0,1[$, thus, 
$S_{(\theta)}(E)(v)\le \theta S_{(1)}(E)(v)$ for every $v\in\rvon$
and every $\theta\in]0,1[$.
Hence, $\Lb_{(\theta)}(E)(u)\llim\limits_{\theta\to 0^+} 0$. On the other hand, by
Lemma 3.3 and (3.1) we obtain 
$\Lb_{(\theta)}(E)(u)\llim\limits_{\theta\to +\infty} +\infty$.
Moreover the map $\theta\mapsto \Lb_{(\theta)}(E)(u)$ is continuous by Lemma 3.1,
and strictly increasing.
The latter property holds as the map $\theta\mapsto \S_{(\theta)}(E)(v)$ is strictly
increasing for every nonconstant $v\in\rvon$. 
Now, the existence and uniqueness of $\bar\theta_{\sigma;E}(u)$ easily follows, 
as well as the continuity of $\bar\theta_{\sigma;E}$. \enpr

\bigskip
We  define   $\ca_2$ to be the set of 
$E\in\ca_1$ such there exists
$\witi E$ 
having  the following properties:

\qua a) $\witi E$ satisfies all properties in the definition of $\ca_1$ except possibly
for the fact that in $Q_4$) we do not require that the inequality is strict  if
$ (u\wedge a)\vee b\ne u$.

\qua b) $\witi E$ is $p$-homogeneous,

\qua c) $\witi E$ is an eigenform 
for $\Lb_{(1)}$, in other words, there exists $\rho>0$ such that
$\Lb_{(1)}(\witi E)=\rho \witi E$,

\qua d) ${E(u)\over\witi E(u)}\llim\limits_{u\to 0}1$.

\bigskip
{\bf Lemma 3.6.} {\sl If $E\in\ca_2$ and $c$ is a constant function in $\rvo$,
 then $\bar\theta_{\sigma,E}(u)\llim\limits_{u\to c} 
\disp{\big({\sigma\over\rho}\big)^{1\over p}}$.}
\smad
Given $\veps>0$ let
$\delta>0$ be such that

$${1\over 1+\veps}< {E(u)\over \witi E(u)}< 1+\veps \ \text{if}\ u\in\rvo, u
\  \text{ nonconstant\ and}
\max|u|<\delta. \eqno(3.3)$$

Of course (3.3) holds also if $u\in\rvo$ is  nonconstant and $\osc(u)<\delta$, as in such  a case,
in view of $Q_2$),
we can replace $u$ by $u-u(P_j)$.
Take now $u$ 
with $\disp{\osc(u)<\min\big\{\delta,{\delta\over\theta}\big\}}$,
$\theta=\disp{\big({t\over \rho}\big)^{1\over p}}$, $t>0$. We have

$$\Lb_{(\theta)}(\witi E)(u)
={t\over\rho} \Lb_{(1)}(\witi E)(u)=t \witi E(u).$$

Let $v\in H'_{\witi E,(\theta)}(u)$ and use $t=\sigma(1+\veps)^{-2}$.
We have $\osc(\theta v)\le\theta\osc(u)<\delta$, and 

$$\Lb_{(\theta)}(E)(u)\le S_{(\theta)}(E)(v)\le (1+\veps) S_{(\theta)}(\witi E)(v)=
(1+\veps) \Lb_{(\theta)}(\witi E)(u)$$
$$=(1+\veps)t \witi E(u)\le t(1+\veps)^2  E(u)=\sigma E(u)$$

where in the inequalities comparing $E$ and $\witi E$ we have used (3.3),
so that $\bar\theta_{\sigma,E}(u)\ge \disp{\big({\sigma\over \rho(1+\veps)^2}\big)^{1\over p}}$.
Let now $v\in H'_{E,(\theta)}(u)$ and use $t=\sigma(1+\veps)^2$. 
We have $\osc(\theta v)\le\theta\osc(u)<\delta$, and 

$$\Lb_{(\theta)}(E)(u)= S_{(\theta)}(E)(v)\ge (1+\veps)^{-1} S_{(\theta)}(\witi E)(v)\ge
(1+\veps)^{-1} \Lb_{(\theta)}(\witi E)(u)$$
$$=(1+\veps)^{-1} t \witi E(u)\ge t  (1+\veps)^{-2} E(u)=\sigma E(u)$$

so that $\bar\theta_{\sigma,E}(u)\le\disp{\big({\sigma (1+\veps)^2\over \rho}\big)^{1\over p}}$.
In conclusion, $\bar\theta_{\sigma,E}(u)\llim\limits_{u\to c} 
\disp{\big({\sigma\over\rho}\big)^{1\over p}}$. \enpr

\bigskip
If $E\in \ca_{2}$, we now put $\bar\theta_{\sigma,E}(u)=
\disp{\big({\sigma\over\rho}\big)^{1\over p}}$ if $u$ is constant, so that, by 
Lemma 3.3
 the map $\bar\theta_{\sigma,E}$ is continuous on
all of $\rvo$. Also, if $E\in\ca_2$ (or more generally
if $E\in\ca_1$ provided $u$ is nonconstant) define

$$H_{\sigma;E}(u)=H_{E,\big(\bar\theta_{\sigma,E}(u)\big)}(u)$$

 and similarly 
 $H'_{\sigma;E}(u)=H'_{E,\big(\bar\theta_{\sigma,E}(u)\big)}(u)=
H_{\sigma,E}(u)\cap \cl'(u)$,
  and $\bar S(E)(v)=S_{\big(\bar\theta_{\sigma,E}(v)\big)}(v)$.
Moreover, let $\bar H'_{\sigma;E}(u)$ be a specific element of
$H'_{\sigma;E}(u)$ (which we choose arbitrarily).

We now define a new subclass  of $\ca_1$.
Namely, let $\ca_3$ be
 the set of $E\in\witi\ca_1$ such that

\qua $Q_5$) If $u,v\in\rvo$, and $v(P)\ge 0$ for $P$ such that $u(P)=\min u$
and $v(P)= 0$ for $P$ such that $u(P)>\min u$, then ${d\over dt}^+|_{t=0} E(u+tv)\le 0$.

\smallskip
We put $\ca_4=\ca_2\cap\ca_3$. The following Lemma is in some sense a  variant of $Q_5$), where
we assume stronger conditions on $u$ and $v$ and in fact  does not depend on $Q_5$).
%\qua $P_6$) There exist $\delta>0$ and  $\bar\alpha:]0,\delta[\to\bre$ 
%such that ${\bar\alpha(t)\over t}\llim\limits_{t\to 0} 0$ and 
%$E(u)\le\bar\alpha\big(\osc(u)\big)$.

\bigskip
{\bf Lemma 3.7.} {\sl If $E\in\ca_1$ and 
$u,v\in\rvo$, $u$ nonconstant, and for some $\theta>0$ we have

\quad  $v(P)=1$ for $P$ such that $u(P)=\min u$ and $v(P)=0$ for the other $P$,

then ${d\over dt}^+|_{t=0} E(u+tv)< 0$.}
\smad We can and do assume $\theta=1$, as the general case can be easily
reduced to this one.
Let $m=:\min u$, $m':=\min \{u(P_j):u(P_j)>m\}$, and take
  $\bar t\in ]0, m'-m[$. Let 
$A=\{P\in\vo:U(P)=\min u\}$. Then, we have
$u+\bar t v=(u\wedge\max u)\vee(\min u+\bar t)$, thus, by $Q_4$) 
we have $E(u+\bar tv)< E(u)$. As $E$ is convex, then

$$\disp{{d\over dt}^+|_{t=0} E(u+tv)\le {E(u+\bar t v)-E(u)\over \bar t}<0}. \enpr $$

\bigskip
{\bf Lemma 3.8.} {\sl If $E\in\ca_3$ and $u$ is nonconstant, then  
for every $\theta>0$ we have $H_{E,(\theta)}(u)\subseteq\cl''(u)$.}
\smad If $v\in\cl(u)\setminus\cl'(u)$, then 
$\tilde v:= (v\wedge\max u)\vee\min u\in \cl'(u)$ and, by $Q_4$) we have
$S_{(\theta)}(E)(\tilde v)<S_{(\theta)}(E)(v)$. Hence, it suffices to prove that, 
for every given $v\in\cl'''(u)$, then 

$$\exists\, w\in\cl(u):
S_{(\theta)}(E)(w)<S_{(\theta)}(E)(v).\eqno (3.4)$$

Let

$$w_t:=v+t\chi_B,\quad  B:=\{Q\in\von\setminus\vo: v(Q)=\min u\},$$

for positive $t$. Then,

$$w_t\circ\psi_i=v\circ\psi_i + t s_i,\quad s_i(P_j):=\cases
1&\ \text{if}\ j\ne i, v\big(\psi_i(P_j)\big)=\min u\\ 0&\ \text{otherwise.}\endcases$$

By $Q_5$) with $\theta v\circ\psi_i$ in place of $u$ and
$\theta s_i$ in place of $v$, we have ${d\over dt}^+|_{t=0} E(\theta w_t\circ \psi_i)\le  0$ for all $i$. 
On the other hand, by the definition of $\cl'''(u)$, 
there exists $\bar i$ such that $v\circ \psi_{\bar i}$ attains both 
$\max u$ and $\min u$, so that $\min  v\circ \psi_{\bar i}=\min u$
and $\max  v\circ \psi_{\bar i}=\max u$. Moreover, clearly, at least one of a) or b) holds, where

\qua a) $\theta  v\circ\psi_{\bar i}$ does not attain its minimum at a point of the form $P_{\bar i}$,

\qua b) $\theta  v\circ\psi_{\bar i}$ does not attain its maximum at a point of the form $P_{\bar i}$.

Of course, if $\bar i>N$, the point $P_{\bar i}$ does not exist, and in such a case both a) and b) hold.
We are going to prove, that if a) holds, then also (3.4) holds. The case b) can be
 reduced to a) replacing $u$ by $-u$ and $v$ by $-v$. By a),
in view of Lemma 3.7 with $\theta v\circ\psi_{\bar i}$ in place of $u$ and
$\theta s_{\bar i}$ in place of $v$, we have ${d\over dt}^+|_{t=0} E(\theta 
w_t\circ \psi_{\bar i})< 0$.  In conclusion, 
${d\over dt}^+|_{t=0} S_{(\theta)}(E)(w_t)< 0$, 
thus for some positive $t$ we have $S_{(\theta)}(E)(w_t)<S_{(\theta)}(E)(v)$. 
Moreover, $w_t\in\cl(u)$ as $\chi_B=0$ on $\vo$.
\enpr

\bigskip

\centerline{\bf 4. Energies on the fractal.}

\medskip

In Section 2,
we saw that  the functionals $\witi S_n(E)$ increasingly converge
to a self-similar energy defined on all of the fractal when $E\in\td$. 
In this Section, we will  study 
analogous notions when $E\in\ca_1$, but  based on $S_{(\theta)}$.
Recall that, for every $\theta>0$ we have defined

$$S_{(\theta)}(E)(v)=\sum\limits_{i=1}^k E(\theta v\circ\psi_i)$$

 when
$E:\rvo\to[0,+\infty[$, and, when $\bte\in\Theta$, $\Theta$ denoting the set 
of the continuous functions from $\rvo$ into $]0,+\infty[$, 
 we can also define

$$\witi S_{\sigma,\bte}(E)(v)=
{1\over\sigma} S_{(\bte(v|_{\vo}))}(E)(v),$$

 if $v\in\rvon.$
In such a case $S_{(\theta)}(E)$ is defined from $\rvon$ to $\bre$.
More generally, consider $\ce:\rvin n\to\bre$ $n=0,1,2,,,\infty,\omega$ with the
comventions $\vin {\omega}=\cf$, $\infty+n=\infty$, $\omega+n=\omega$.
Then if $\sigma>0$ and $\bte\in\Theta$, 
 we define $S_{(\theta)}(\ce):\rvin {n+1}\to\bre$   by

$$S_{(\theta)}(\ce)(v)=\sum\limits_{i=1}^k \ce(\theta v\circ\psi_i),
\quad \witi S_{\sigma,\bte}(\ce)(v)=
{1\over\sigma} S_{(\bte(v|_{\vo}))}(\ce)(v)$$

for every $v\in \rvin {n+1}$.  If $\bte\in\Theta$,
 define now $\bte_{,i_1,...,i_n}(v)$ 
 for $v:\vin n\to\bre$ by recursion as
 
 $$\bte_{,}(v)=1,\ \  \bte_{,i}(v)=\bte(v|_{\vo}),$$
$$\bte_{,i_1,...,i_n,i_{n+1}}(v)=\bte_{,i_2,...,i_{n+1}}\big(\bte(v|_{\vo}) 
 v\circ\psi_{i_1}\big) \bte(v|_{\vo}).$$

When $m=0,1,2,3,....,\infty,\omega$, $\ce:\rvin m\to [0,+\infty[$ and  $\bte\in\Theta$, we now define 
$S_{n;\bte}(\ce)$ from $\rvin {n+m}$ into $[0,+\infty[$ as follows

$$S_{n;\bte}(\ce)(v)=\sum\limits_{i_1,...,i_n=1}^k \ce  (\bte_{,i_1,...,i_n}(v) \, v
\circ \psi_{i_1,...,i_n})\quad \forall  v\in\rvin {n+m}, $$
$$\witi S_{n;\sigma,\bte}(\ce)={1\over\sigma^n} S_{n;\bte}(\ce),
\quad \witi\ce_{n,\sigma,E}=\witi S_{n;\sigma,\bte_{\sigma,E}}(E)$$

%$$\witi \ce_{1,\sigma,E}(v)=\witi S_{\sigma,E}(\ce_E)(v).\eqno (4.2)$$

Note that $S_{0;\bte}(\ce)=\ce$. The following property will be useful.

\bigskip
{\bf Lemma 4.1} {\sl We have
$\witi\ce_{1,\sigma,E}(v)\ge E(v|_{\vo})$ for every $v\in\rvon$ and the equality
holds if and only if $v\in H_{\sigma,E}(v|_{\vo})$.}
\smad 
By definition,  for every $\theta>0$ 
we have $S_{(\theta)}(E)(v)\ge \Lb_{(\theta)}(E)(v|_{\vo})$ and the equality holds if and only if
$v\in H_{E,(\theta)}(v|_{\vo})$. To obtain the Lemma it suffices to
use $\theta=\bte_{\sigma;E}$. \enpr

\bigskip
{\bf Lemma 4.2.} {\sl Given $\ce$ and $\bte$  as above, we have 

\qua i) If $\ce$ is continuous, then $\witi S_{n;\sigma,\bte}(\ce)$ is  continuous on $\rvin {n+m}$.

\qua ii)  $\witi S_{n+1;\sigma,\bte}(\ce)=  
\witi S_{\sigma,\bte}\big(\witi S_{n;\sigma,\bte}(\ce)\big) $ on $\rvin {n+m+1}$.

\qua iii)   $\witi S_{n;\sigma,\bte}(\ce)=  \witi S_{\sigma,\bte}^n (\ce) $
on  $\rvin {n+m}$.

\qua iv) We have $\witi \ce_{n+1,\sigma, E}(v)\ge \witi \ce_{n,\sigma, E}(v)$  for every $v\in  \rvin {n+1}$.}

\smad  i) follows from  the continuity of
$\bte$. We prove ii). In view of (4.1), (4.2)
 and the definition of $\bte_{\sigma,E}$, we have 

$$\witi S_{n+1;\sigma,\bte}(\ce)(v)={1\over\sigma^{n+1}}\sum\limits_{i_1,...,i_n,
i_{n+1}=1}^k \ce (\bte_{,i_1,...,i_n,i_{n+1}}(v) \, v
\circ \psi_{i_1,...,i_n, i_{n+1}})$$
$$={1\over\sigma}\ \sum\limits_{i_{1}=1}^k
{1\over\sigma^n}\sum\limits_{i_2,...,i_{n+1}=1}^k \, 
\ce\Big(\bte_{,i_2,...,i_n,i_{n+1}}\big(\bte(v|_{\vo})\,  v\circ\psi_{i_1}\big)
\bte(v|_{\vo}) v \circ \psi_{i_1}\circ\psi_{i_2,...,i_{n+1}}\Big)=$$
$$= {1\over\sigma}\sum\limits_{i_{1}=1}^k
 \witi S_{n;\sigma,\bte}(\ce)(\bte(v|_{\vo}) v\circ\psi_{i_1})=
 \witi S_{\sigma,\bte}\big(\witi S_{n;\sigma,\bte}(\ce)\big)(v)$$

We have so proved ii) and iii) follows immediately from ii).
iv) By Lemma 4.1,  iv) holds for $n=0$. The general case, in view of ii), follows by recursion as
the map $\ce\mapsto \witi S_{\sigma,\bte}(\ce)$ is increasing.
\enpr

\bigskip
Put now $\witi \ce_{\infty,\sigma,E}(v):=\lim\limits_{n\to +\infty}
 \witi \ce_{n,\sigma,E}(v)$ for every $v\in\rvin {\infty}$. It easily follows
from Lemma 4.2 iv) that such a limit exists. We equip  $\rvin {\omega}$
with the norm $||\ ||_{\infty}$ defined as usual as $||v||_{\infty}=\sup |v|$. 
We could do the same on  $\rvin {\infty}$. However, strictly speaking,
this is not a norm on $\rvin {\infty}$ in that it can attain the value $+\infty$.

\bigskip
{\bf Lemma 4.3.} {\sl The functional $\witi \ce_{\infty,\sigma,E}$ on  $\rvin {\omega}$ with the norm
$||\ ||_{\infty}$ is lower semicontinuous. Moreover,  it  is self-similar in the sense that

$$\witi \ce_{\infty,\sigma,E}=  \witi S_{\sigma,\bte_{\sigma,E}}
\big( \witi \ce_{\infty,\sigma,E}\big). $$}
\smad 
The map $v\mapsto v|_{\rvin n}$ is continuous from
$\rvin {\omega}$ into $\rvin n$, hence by Lemma 4.2 i) $\witi \ce_{\infty,\sigma,E}$
is the supremum of continuous functionals, thus is lower semicontinuous. The selfsimilarity follows from
Lemma 4.2 iv. In fact, for every $v\in\rvin {\omega}$ we have

$$ \witi S_{\sigma,\bte_{\sigma,E}}
\big( \witi \ce_{\infty,\sigma,E}\big)(v)={1\over\sigma}\sum\limits_{i=1}^k
\witi \ce_{\infty,\sigma,E}\big(\bte_{\sigma,E}(v|_{\vo}) v\circ\psi_i\big)$$ 
$$=\lim\limits_{n\to +\infty}
{1\over\sigma}\sum\limits_{i=1}^k
\witi \ce_{n,\sigma,E}\big(\bte_{\sigma,E}(v|_{\vo}) v\circ\psi_i\big)
=\lim\limits_{n\to +\infty}\witi S_{\sigma,\bte_{\sigma,E}}
\big( \witi \ce_{n,\sigma,E}\big)(v)$$
$$=\lim\limits_{n\to +\infty} \witi \ce_{n+1,\sigma,E}(v)=\witi \ce_{\infty,\sigma,E}(v). \enpr
$$

\bigskip

\bigskip
\centerline{\bf 5. Some Lemmas.}

\medskip
{\bf Lemma 5.1.} {\sl Given $E\in\ca_3$, for every real numbers $a, b$ with $0<a\le b$, 
there exists $\alpha_{a,b}\in ]0,1[$ such that, if $u\in\rvo$ and $\osc (u),\theta\in [a,b]$,
then every $v\in H_{E,(\theta)}(u)$ satisfies

$$\osc(v\circ\psi_i)\le \alpha_{a,b}\, \osc(u)\quad \forall\, i=1,...,k.$$}
 \smad
 By contradiction, if the Lemma is false, we find $u_n\in\rvo$ and we can and do assume
 $u_n(P_1)=0$, and $\theta_n$ such that $\osc (u_n),\theta_n\in [a,b]$
 and $v_n\in H_{E,(\theta_n)}(u_n)$ and $i_n=1,..,k$ such that
 $\disp{ {\osc(v_n\circ\psi_{i_n})\over\osc(u_n)}\llim\limits_{n\to+\infty} 1}$. 
 By taking a subsequence we can and do assume that $\theta_n\llim\limits_{n\to+\infty} 
 \theta\in [a,b]$,  $u_n\llim\limits_{n\to+\infty}  u$,
  $v_n\llim\limits_{n\to+\infty} v$ and $i_n=i$. By continuity,
 we have  $\osc(u)\in [a,b]$,  so that $u$ is nonconstant, $v\in H_{E,(\theta)}(u)$ 
 and $\osc(v\circ\psi_i)=\osc(u)$. This contradicts Lemma 3.8. \enpr

\bigskip
{\bf Lemma 5.2.} {\sl Let $E\in\ca_2$
and let $\witi E$ be as in the definition of $\ca_2$.Then, for every $\sigma>0$
 there exist $\bar\alpha\in ]0,1[$ and $\bar\delta>0$ such that   
if $u\in\rvo$   with $\osc(u)< \bar\delta$, and
$v \in H'_{\sigma,E}(u)$,  we have

$$ \witi E(v\circ\psi_i)\le \bar\alpha  \witi E(u)\quad \forall\, i=1,...,k..$$}
\smad
Given $\eta>0$,
there exists $\delta_{\eta}>0$ such that if $\osc(u)<\delta_{\eta}$, then

$$\big({1\over 1+\eta}\big)\witi E(u)\le E(u)\le (1+\eta)\witi E(u).$$

By Lemma  3.6 we deduce that
there exist $\delta'>0$ such that, if $\osc(u)< \delta'$, then, 
  as $\rho<1$,
we have $\bar\theta_{\sigma,E}(u)\in [\theta_1,\theta_2]$
 with $\theta_1>\sigma^{1\over p}$.  We now choose
$\eta>0$ so small that $(1+\eta)^2<{\theta_1^p\over\sigma}$ and
$\bar\delta=\min\big\{\delta', \delta_{\eta},\disp{\delta_{\eta}\over\theta_2}\big\}$.
 Suppose  $u\in\rvo$ 
 with $\osc(u)<\bar\delta$, and
$v\in H'_{\sigma,E}(u)$, and let $\theta=\bar\theta_{\sigma,E}(u)$. 
Note that   $\osc(u)<\delta_{\eta}$, hence $E(u)\le (1+\eta)\witi E(u)$.
Moreover, for every $i=1,...,k$ we have

$$\osc(\theta v\circ\psi_i)=\theta \osc(v\circ\psi_i)\le\theta \osc(u)\le\theta_2\osc(u)<\delta_{\eta}.$$
Thus

$$\theta^p\witi E(v\circ\psi_i)
=\witi E(\theta v\circ\psi_i)$$
$$\le  (1+\eta)  E(\theta v\circ\psi_i)
\le (1+\eta)  S_{(\theta)}(E)(v)= $$
$$\sigma(1+\eta)E(u)\le \sigma(1+\eta)^2 \witi E(u),$$

so that putting $\bar\alpha:= \disp{ \sigma (1+\eta)^2\over \theta_1^p} \in ]0,1[$, we have

$$\witi  E(v\circ\psi_i)\le \bar\alpha  \witi E(u).\enpr$$

\bigskip
{\bf Lemma 5.3.} {\sl If $\sigma\in ]0,1]$, $E\in\ca_2$, $v\in\rvin {\infty}$ and
$\witi\ce_{\infty,\sigma,E}(v)<+\infty$, then
there exists $H_v>0$ such that

$$\osc\Big(\bte_{,i_1,...,i_n}(v) \, v\circ \psi_{i_1,...,i_n}|_{\vo}\Big)\le H_v$$

where $\bte=\bte_{\sigma,E}$, for all $n$ and $i_1,...,i_n$.}
\smad
We have

$$E (\bte_{,;i_1,...,i_n}(v) \, v\circ \psi_{i_1,...,i_n})\le 
\sigma^n\witi\ce_{n,\sigma,E}(v)<\witi\ce_{\infty,\sigma,E}(v)<+\infty$$

and, in view of Lemma 3.3,  we conclude. \enpr
\bigskip

\centerline{\bf 6. The main results.}

\medskip

{\bf Theorem 6.1.} 
{\sl If $E\in\ca_4$, then for every 
$\sigma\in ]0,1]$ and $u\in\rvo$ there exists $v\in\cl'(\cf,u)$ 
such that  $\witi\ce_{\infty,\sigma,E} (v)=E(u)$.}
\smad
Let $v_0:=u$. We are going to  prove  that, given   $v_n\in\rvin n$  satisfying

$$\witi\ce_{n,\sigma,E} (v_n)=E(u), \quad \min u\le v_n(Q)\le\max u\ \forall\, Q\in\vin n,\eqno (6.1)$$

we can extend $v_n$ to a function $v_{n+1}\in \rvin {n+1}$
satisfying (6.1) with $n+1$ in place of $n$.
Put $\bte:=\bte_{\sigma,E}$ and, for every $i_1,...,i_n=1,...,k$, we define 

$$w_{i_1,...,i_n}=\bar H'_{\sigma,E}\big( \bte_{,i_1,...,i_n}(v_n) v_n\circ\psi_{i_1,...,i_n}\big)
$$

and $v_{n+1}\in\rvin {n+1}$ in the following way. If $Q\in\vin {n+1}$, then we have
$Q=\psi_{i_1,...,i_n}(Q_1)$ for some $Q_1\in\von$ and $i_1,...,i_n=1,...,k$, and we put

$$v_{n+1}(Q)={w_{i_1,...,i_n}(Q_1)\over  \bte_{,i_1,...,i_n}(v_n)}.$$

Such a definition is correct. We have to prove that, if $Q=\psi_{i_1,...,i_n}(Q_1)=\psi_{i'_1,...,i'_n}(Q'_1)$,
then $\disp{{w_{i_1,...,i_n}(Q_1)\over  \bte_{,i_1,...,i_n}(v_n)}=
{w_{i'_1,...,i'_n}(Q'_1)\over  \bte_{,i'_1,...,i'_n}(v_n)}}$. To prove this, note that by $P_4$)
$Q_1,Q'_1\in\vo$, thus by the definition of $w_{i_1,...,i_n}$ we have

$$\disp{{w_{i_1,...,i_n}(Q_1)\over  \bte_{,i_1,...,i_n}(v_n)}=
v_n\circ\psi_{i_1,...,i_n}(Q_1)=v_n(Q)=}$$
$$\disp{v_n\circ\psi_{i'_1,...,i'_n}(Q'_1)={w_{i'_1,...,i'_n}(Q'_1)\over  \bte_{,i'_1,...,i'_n}(v_n)}}.
$$

Also, if $Q\in\vin n$, then we can choose $Q_1\in \vo$, and the above argument shows that $v_{n+1}(Q)=v_n(Q)$, so that $v_{n+1}$ is in fact an extension of $v_n$. Next, 

$$\bte_{,i_1,...,i_n}(v_n) v_{n+1}\circ\psi_{i_1,...,i_n}=\bar H'_{\sigma,E}
\big( \bte_{,i_1,...,i_n}(v_n) v_n\circ\psi_{i_1,...,i_n}\big). \eqno (6.2)$$

In particular, for every $Q\in\von$ we have

$$\min\limits_{\vo} v_n\circ \psi_{i_1,...,i_n}\le v_{n+1}\circ \psi_{i_1,...,i_n}(Q)\le 
\max\limits_{\vo} v_n\circ \psi_{i_1,...,i_n}\eqno (6.3)$$

Next, in view of Lemma 4.1, (6.2) and Lemma 4.2 ii), we get

$$\witi\ce_{n+1,\sigma;E}(v_{n+1})=\witi S_{n,\sigma,E} (\witi\ce_{1,\sigma;E})(v_{n+1})=$$
$${1\over \sigma^n}\sum\limits_{i_1,...,i_n=1}^k 
\witi\ce_{1,\sigma;E}\big(\bte_{,i_1,...,i_n}(v_n) v_{n+1}\circ\psi_{i_1,...,i_n}\big)=$$
$${1\over \sigma^n}\sum\limits_{i_1,...,i_n=1}^k 
E\big(\bte_{,i_1,...,i_n}(v_n) v_n\circ\psi_{i_1,...,i_n}\big)=$$
$$\witi S_{n,\sigma,E} (E)(v_{n})=\witi\ce_{n,\sigma;E}(v_{n}).$$

Thus, taking also into account (6.3), $v_{n+1}$ satisfies (6.1) with $n+1$ in place of $n$. As
(6.1) is trivially satisfied for $n=0$, we have constructed a function $v$ on $\vin {\infty}$ as the extension
of every $v_n$, and, in order to prove the Theorem, it suffices to show that $v$ that be extended continuously
on $\cf$, and in order to obtain this, it suffices in turn to prove that
$v$   is uniformly continuous on $\vin{\infty}$.
Hence, in view of Lemma 2.1 i)
it suffices to prove that 

$$\max\big\{\osc(v\circ\psi_{i_1,...,i_n}): i_1,...,i_n=1,...,k\big\}
\llim\limits_{n\to +\infty} 0. \eqno (6.4)$$

Here, of course, the oscillation is on $\vin {\infty}$, but in view of (6.3), we have

$$\osc(v\circ\psi_{i_1,...,i_n})=\osc_{\vo}(v\circ\psi_{i_1,...,i_n})\quad\forall\, i_1,...,i_n=1,...,k$$ 

and also,
 
 $$\osc(v\circ\psi_{i_1,...,i_n,i_{n+1}})\le \osc(v\circ\psi_{i_1,...,i_n})\quad\forall\, i_1,...,i_n,
i_{n+1}=1,...,k.\eqno (6.5)
 $$

 We are going to prove (6.4). If $\osc(u)=0$, then,  in view of (6.5), $\osc(v\circ\psi_{i_1,...,i_n})=0$
for every $i_1,...,i_n=1,...,k$ and (6.4) is trivial. Thus, suppose $\osc(u)>0$, 
 fix $\veps\in ]0,\osc(u)]$ and consider a  sequence $(i_1,i_2,...,i_n)$. Let 

$$F_{i_1,...,i_n}:=\Big\{h< n: \osc\big(\bte_{,i_1,...,i_h}(v) \, v\circ 
\psi_{i_1,...,i_h}\big)\ge \bar\delta\Big\}.$$

where $\bar\delta$ is defined in Lemma 5.2.
Then, by  Lemma 5.3, (6.2)  and Lemma 5.1, if $h\in F_{i_1,...,i_n}$ we have

$$\osc\big(\theta_{\sigma;i_1,...,i_h}(v) \, v\circ \psi_{i_1,...,i_h,i_{h+1}} \big)\le
\alpha\, \osc\big(\theta_{\sigma;i_1,...,i_h}(v) \, v\circ \psi_{i_1,...,i_h} \big)
$$

for suitable $\alpha=\alpha_{a,b}$ where
$a=\min\{\bar\theta_1,\bar\delta\}$, $b=\max\{\bar\theta_2, H_v\}$,

$$\bar\theta_1:=\min\big\{\bar\theta_{\sigma,E}(w): \osc(w)\in [\bar\delta,H_v]\big\},
\quad\bar\theta_2:=\max\big\{\bar\theta_{\sigma,E}(w): \osc(w)\in [\bar\delta,H_v]\big\}.$$

Therefore, if $h\in F_{i_1,...,i_n}$ we have

$$\osc\big(v\circ \psi_{i_1,...,i_h,i_{h+1}} \big)\le
\alpha\, \osc\big(v\circ \psi_{i_1,...,i_h} \big).
$$

Thus, if  $F_{i_1,...,i_n}$ has at least 

$$C:= \big[\disp{\log_{{1\over \alpha}}{\osc(u) \over\veps}}\big]+1$$

 elements, in view also of (6.5),
we have $\osc(v\circ\psi_{i_1,...,i_n})<\veps$.
On the other hand, if we have 
$s$  consecutive elements $h,...,h+s-1$  in $\{0,1,...,n-1\}\setminus F_{i_1,...,i_n}$
by Lemma 5.2, (6,2) and the homogeneity of $\witi E$ we have

$$\witi E(v\circ\psi_{i_1,...,i_{h+s}})\le \bar\alpha^s
\witi E(v\circ\psi_{i_1,...,i_{h}}).\eqno (6.6)$$

Now, let $m,M$  be defined as

$$m=\min\big\{\witi E(w):\osc(w)\ge\veps\big\},\quad
M=\max\big\{\witi E(w):\osc(w)\le \osc(u)\big\}.
$$

If $\bar s$ is such that ${m\over M}>\bar\alpha^{\bar s}$,  as
$\osc\big(v\circ\psi_{i_1,...,i_{h}}\big)\le \osc(u)$, then
$\witi E\big(v\circ\psi_{i_1,...,i_{h}}\big)\le M$, hence
by (6.6)  $\witi E\big(v\circ\psi_{i_1,...,i_{h+\bar s}}\big)<m$ and 
$\osc \big(v\circ\psi_{i_1,...,i_{h+\bar s}}\big)<\veps$, so that, by (6.5)
$\osc \big(v\circ\psi_{i_1,...,i_n}\big)<\veps$.
If $n>C\bar s$, then
for every $i_1,...,i_n=1,...,k$,  either $F_{i_1,...,i_n}$ has at least $C$ elements,
or there exist $\bar s$ consecutive elements
in $\{0,1,...,n-1\}\setminus F_{i_1,...,i_n}$, thus in any case
 $\osc(v\circ \psi_{i_1,...,i_n})<\veps$, and
we have proved (6.4).  \enpr

\bigskip
{\bf Theorem 6.2.}
{\sl  For every $E\in\ca_4$ and every $\sigma\in]0,1]$,
 $\witi\ce_{\infty,\sigma,E}$ is the only functional $\ce$  from
$\bre^{\cf}$ to $\bre$ satisfying the following
\smallskip

\qua a) $\ce$ is lower semicontinuous on $\bre^{\cf}$
 with respect to the $L_{\infty}$ metric

\qua b) there exists $\bte\in\Theta$  such that
$\ce=\witi S_{\sigma,\bte}(\ce)$.

\qua c) For every $u\in\rvo$ there exists $v\in\cl'(\cf,u)$ 
 such that $E(u)=\ce(v)=
\min\limits_{w\in \cl(\cf,u)} \ce(w)$.}
\smad
By Lemma 4.3 $\witi\ce_{\infty,\sigma,E}$ satisfies a) and  b).
Moreover, by Lemma 4.2 iv) and Theorem 6.1, it also satisfies c).
I will now prove that $\witi\ce_{\infty,\sigma,E}$ is the unique functional 
satisfying a), b) and c). Suppose $\ce$ satisfies a), b) and c).
Note that $\theta_{,i_1,...,i_n}(v)$ only 
depends on $v|_{\vin {n-1}}$, and this can be proved by recursion.
Now, as b) holds, then

$$\ce(v)=\witi S_{\sigma,\bte}^n(\ce)(v)=\witi S_{n;\sigma,\bte}(\ce)=
{1\over\sigma^n}\sum\limits_{i_1,...,i_n=1}^k
 \ce\big(\bte_{,i_1,...,i_n}(v) v\circ\psi_{i_1,...,i_n}\big) \eqno(6.7)$$

for every $v\in\bre^{\cf}$, every $n$ and every $i_1,...,i_n=1,...,k$. 
Given $v\in\bre^{\cf}$,  for every $i_1,...,i_n$ let $v_{i_1,...,i_n}\in\cl'(\cf, v\circ\psi_{i_1,...,i_n})$
 be  such that

$$E\big(\bte_{,i_1,...,i_n}(v) v|_{\vin n}\circ\psi_{i_1,...,i_n}\big)=
\ce\big(\theta_{,i_1,...,i_n}(v) v_{i_1,...,i_n}\big)\le E(w)   \eqno (6.8)$$
for every  $w\in \bigl(\cf, 
\bte_{,i_1,...,i_n}(v) v|_{\vin n}\circ\psi_{i_1,...,i_n}\big)$.
Moreover, define $v_n$ on $\cf$ by

$$v_n\circ\psi_{i_1,...,i_n}=v_{i_1,...,i_n}\quad\forall\, i_1,...,i_n=1,...,k$$

By an argument similar to that used in the proof of Theorem 6.1 we see that such a definition is correct
and that $v_n=v$ on $\vin n$, hence $\bte_{,i_1,...,i_n}(v)=\bte_{,i_1,...,i_n}(v_n)$.
It follows that
$\ce(v)\ge\ce(v_n)$  for every $ n$. In fcat, using (6.7) both for $v$ and for $v_n$ and (6.8) we get

$$\ce(v)={1\over\sigma^n}\sum\limits_{i_1,...,i_n=1}^k
 \ce\big(\bte_{,i_1,...,i_n}(v) v\circ\psi_{i_1,...,i_n}\big)\ge$$
$${1\over\sigma^n}\sum\limits_{i_1,...,i_n=1}^k
 \ce\big(\bte_{,i_1,...,i_n}(v_n) v_n\circ\psi_{i_1,...,i_n}\big)=\ce(v_n).$$

Hence,
$\ce(v)\ge\liminf\limits_{n\to +\infty}\ce(v_n).$
On the other hand,
as for every $Q\in \cf_{i_1,...,i_n}$ we have
$v(Q), v_n(Q)\in  [\min\limits_{\cf_{i_1,...,i_n}} v, \max\limits_{\cf_{i_1,...,i_n}} v]$,
then, by Lemma 2.1 ii),
 $v_n\llim\limits_{n\to +\infty} v$ uniformly on $\vin {\infty}$, but by a continuity argument also
on $\cf$, hence by a)
$\ce(v)\le\liminf\limits_{n\to +\infty}\ce(v_n).$ Thus

$$\ce(v)=\liminf\limits_{n\to +\infty}\ce(v_n).$$

On the other hand, we have
$$\ce(v_n)={1\over\sigma^n}\sum\limits_{i_1,...,i_n=1}^k \ce
\big(\bte_{,i_1,...,i_n}(v) v_n\circ\psi_{i_1,...,i_n}\big)$$
$$={1\over\sigma^n}\sum\limits_{i_1,...,i_n=1}^k E
\big(\bte_{,i_1,...,i_n}(v) v|_{\vin n}\circ\psi_{i_1,...,i_n}\big)$$

so that $\ce(v)$ is determined by $E$. \enpr

\bigskip\bigskip
\centerline {\bf References }

%\bigskip
%\cite{1}\ B.M. Hambly, V. Metz, A. Teplyaev,
%{\it Self-Similar Energies on Post-critically Finite Self-Similar Fractals}, J. London Math. 
%Soc. {\bf 74} (2006) 93-112.

\smallskip
\cite{1}\ K. Hattori, T. Hattori, H. Watanabe,
 {\it Gaussian Field Theories on General Networks and the Spectral 
Dimension},  Progr. Theor. Phys. Suppl. {\bf 92} (1987),
108-143. 

\cite{2} P. E. Herman, R. Peirone, Robert Strichartz,
 {\it $p$-Energy and $p$-Harmonic Functions on Sierpinski Gasket
Type Fractals},
 {\it Potential Analysis},    {\bf 20},  2 (2004),  125-148.

\smallskip
\cite{3}\ J. Kigami, {\it Harmonic Calculus on p.c.f. Self-similar  Sets},
Trans. Amer. Math. Soc. {\bf 335} (1993),   721-755.

\smallskip
\cite{4}\ J. Kigami, {\it Analysis on fractals}, 
Cambridge Tracts in Mathematics, 143.
Cambridge University Press, Cambridge, 2001.

\smallskip
\cite{5}\ T. Lindstr\o m, {\it Brownian Motion on Nested Fractals},  Mem. Amer. Math. Soc.  
No. 420 (1990).

%\smallskip
%\cite{6}\ V.Metz, {\it The short-cut test},
%Journal of Functional Analysis {\bf 220} (2005), 118-156.

%\smallskip
%\cite{6} R. Peirone, {\it Convergence and Uniqueness Problems for  Dirichlet 
%Forms on Fractals}, Boll. U.M.I., (8), 3-B (2000), 431-460.

\smallskip
\cite{6} R. Peirone, {\it Existence of self-similar
energies on finitely ramified fractals}, to appear on
Journal d'Analyse Mathematique.

%\smallskip
%\cite{7}\ R. Peirone, {\it Existence of Eigenforms on Fractals with three Vertices}, Proc. 
%Royal Soc. Edinburgh,
%{\bf 137 A} (2007),  1073-1080.

%\smallskip
%\cite{8}\ R. Peirone, {\it Uniqueness of Eigenforms on Certain  %Fractals}, 
%preprint del Dipartimento di Matematica, Universit\`a %di Roma-Tor Vergata, 26/10/2006.

%\smallskip
%\cite{8} R. Peirone, {\it Existence of Eigenforms on Nicely Separated Fractals}, 
%in Proceedings of Symposia in Pure Mathematics,
%Amer. Math. Soc., Vol. 77, 231-241, 2008. 

%\smallskip
%\cite{9}\ C. Sabot, {\it Existence and Uniqueness of Diffusions 
%on Finitely Ramified  Self-Similar Fractals},
 %Ann. Sci. \'Ecole Norm. Sup. (4) {\bf 30} (1997), no. 5, 605-673.

%\smallskip
%\cite{11}\ C. Sabot,
%{\it Espaces de Dirichlet reli\'es par des points et application %aux 
%diffusions sur les fractals finiment ramifi\'es},
%Potential Anal. 11 (1999), no. 2, 183--212.

%\smallskip
%\cite{10}\
%R. S. Strichartz, {\it Differential equations on fractals: A tutorial}, 
% Princeton University Press, 2006.

\end

\bigskip
{\bf Remark 3.4.}  We have $\witi\theta_{\sigma,E}\in [\sigma,\sigma l_{\von}]$ 
provided $ \sigma\ge 1$. Here
$l_{\von}$ is by definition the minimum number $l$ such that 
every two points in $\vo$ are connected in $\cg_1$ by a path of length $\le l$. 
In fact, for every $j=1,...,N$,   every $t>0$ and $\theta<\sigma$, we have
$t e_j\in\cl(t e_j)$ and  

$$S_{(\theta)}(E)(t e_j)=\theta E(t e_j)<\sigma E(t e_j)$$

so that $\Lb_{(\theta)}(E)(t e_j)<\sigma E(t e_j)$, and $\bar\theta_{\sigma,E}(t e_j)\ge \sigma$.

\qua On the other hand, for every $t>0$ let $A$ be the nonempty proper subset
of $\vo$ that minimizes $E(t\chi_A)$. Then, if $v\in\cl(t\chi_A)$,  
there exists, by definition of $l_{\von}$, at least one $i=1,...,k$ and $j_1,j_2$ 
such that 

$$ v\big(\psi_i(P_{j_1})\big)-v\big(\psi_i(P_{j_2})\big) \ge
{1\over l_{\von}}t,\quad v\big(\psi_i(P_{j})\big)\notin \big]v\big(\psi_i(P_{j_1})\big),
v\big(\psi_i(P_{j_2})\big)\big[\quad\forall\, j.$$

Therefore, using also ii) and iv), for every $\theta>\sigma l_{\von}$

$$S_{(\theta)}(E)(v)\ge E(\theta v\circ\psi_i)> E(\sigma l_{\von} v\circ\psi_i)=$$
$$E\Big(\sigma l_{\von}\big( (v\circ\psi_i)-v(\psi_i(P_{j_1})\big)\big)\Big)\ge
\sigma E(t\chi_B)\ge \sigma E(t\chi_A)$$

where $B=\big\{P_j: v\big(\psi_i(P_{j})\big)\ge v\big(\psi_i(P_{j_2})\big)\big\}$,
as the function $\sigma l_{\von}\big( (v\circ\psi_i)-v(\psi_i(P_{j_1})\big)\big)$ takes values $\ge \sigma t$ on $B$ 
and $\le 0$ on $\vo\setminus B$.
As a consequence $\Lb_{(\theta)}(E)(t \chi_A)>\sigma E(t \chi_A)$ 
and $\bar\theta_{\sigma,E}(t\chi_A)\le \sigma l_{\von}$. \enpr

\bigskip
{\bf Corollary 5.3.} {\sl Let $E\in\ca_4$. Then, for every $\sigma>0$
and every $b>0$
 there exist $\bar\alpha\in ]0,1[$ and $\bar\delta>0$ such that   
if $u\in\rvo$   with $\osc(u)\le  b$, and
$v \in H'_{\sigma,E}(u)$,  putting
$G_E(u)=\osc(u) E(u)$,
we have

$$G_E(v\circ\psi_i)\le \tilde\alpha  G_E(u)\quad \forall\, i=1,...,k..\eqno(5.1)$$}
\smad
If $\osc(u)<\bar\delta$ then (5.1) holds with $\tilde\alpha=\bar\alpha$ by Lemma 5.2.
If $\osc(u)\in [\bar\delta, b]$, then $\theta_{\sigma,E}\in[\theta_1,\theta_2]$
for some positive $\theta_1$ and $\theta_2$. Hence, (5.1) holds
for some $\alpha$ in place of $\tilde\alpha$. Hence
(5.1) holds in any case with $\tilde\alpha=\max\{\alpha,\bar\alpha\}$. \enpr

\bigskip
{\bf Lemma 5.3.} {\sl If $\sigma>0$, $\theta>0$, $v\in\rvin {n+1}$, 
and $\witi\ce_{n,\sigma,E}(v)=\witi\ce_{n+1,\sigma,E}(v)$, then 
$$\theta_{\sigma;i_1,...,i_n}(v) \, v\circ \psi_{i_1,...,i_n}|_{\von}\in H_{\sigma;E}\big(
 \theta_{\sigma;i_1,...,i_n}(v) \, v\circ \psi_{i_1,...,i_n}|_{\vo}\big).$$}
\smad
First, by Lemma 4.2 iii) we have

 $$\witi \ce_{n+1,\sigma,E}(v)=
\witi S_{\sigma,E}^n(\witi\ce_{1,\sigma,E})(v)=
(\witi\ce_{1,\sigma,E})_{n,\sigma,(E)}(v)$$
$$ ={1\over\sigma^n}
\sum\limits_{i_1,...,i_n=1}^k (\witi\ce_{1,\sigma,E}) (\theta_{\sigma;i_1,...,i_n}(v) \, v
\circ \psi_{i_1,...,i_n})$$

Hence, as  
$\witi\ce_{n,\sigma,E}(v)=\witi\ce_{n+1,\sigma,E}(v)$, thus

$$\sum\limits_{i_1,...,i_n=1}^k \witi\ce_{1,\sigma,E} (\theta_{\sigma;i_1,...,i_n}(v) \, v
\circ \psi_{i_1,...,i_n})=\sum\limits_{i_1,...,i_n=1}^k \ce_E
 (\theta_{\sigma;i_1,...,i_n}(v) \, v\circ \psi_{i_1,...,i_n})$$

and,  as $(\witi\ce_{1,\sigma,E})\ge\ce_E$, we have

$$\witi\ce_{1,\sigma,E} (\theta_{\sigma;i_1,...,i_n}(v) \, v\circ \psi_{i_1,...,i_n})=
 \ce_E (\theta_{\sigma;i_1,...,i_n}(v) \, v\circ \psi_{i_1,...,i_n})
$$

and the Lemma follows. \enpr

\quad v) If $u,v\in\rvo$, $u$ nonconstant, and $v(P)\ge 0$ for $P$ such that $u(P)=\min u$,
$v(P)\le 0$ for $P$ such that $u(P)=\max u$, then ${d\over dt}^+|_{t=0} E(u+tv)\le 0$.

\quad vi) There exist $\delta>0$ and 
 $\bar\alpha:]0,\delta[\to\bre$ such that ${\bar\alpha(t)\over t}\llim\limits_{t\to 0} 0$.

\quad vii) There exists $\witi E:\rvo\to [0,+\infty[$ 
and $\rho>0$ such that $\Lb_{(1)}(\witi E)=\rho\witi E$
and ${E(u)\over \witi E(u)}\llim\limits_{u\to 0} 0$.

Denote by $\ca'$ the set of $E\in\ca$ satisfying

\quad viii) $E(u\circ\sigma)=E(u)$ for every permutation $\sigma$ of $\vo$.

\bigskip
{\bf Lemma 1.} {\sl For every $E\in\ca$, $\theta>0$ and $u\in\rvo$ the sets $H_{E,(\theta)}(u)$ 
and $H'_{E,(\theta)}(u)$ are nonempty compact and convex.}
\smad 
We see that $H_{E,(\theta)}(u)$ is nonempty. By convexity and $E(0)=0$
 we have $E(tu)\ge t E(u)$ for every $t>1$.
Hence, if $u$ is nonconstant

$$E(u)=E\big(u-u(P_1)\big)\ge\osc(u) E\Big( {u-u(P_1)\over || u-u(P_1) ||}\Big)
\ge c E(u)\ \ \ \text{if}\ \osc(u)\ge 1$$

where $c:=\min\{ E(u):U(P_1)=0, \osc(u)=1\big\}$. Hence, 
$S_{(\theta)}(E)(v)\llim\limits_{v\to\infty}+\infty$ on $\cl(u)$, so that,
$S_{(\theta)}(E)$ being continuous, it has a minimum on the closed set $\cl(u)$. 
The compactness and the convexity
of $H_{E,(\theta)}(u)$ are trivial, and from this we deduce also t
he compactness and the convexity
of $H'_{E,(\theta)}(u)$ as also $\cl(u)$ is compact and convex. 
Finally, by iv), if $w\in H_{E,(\theta)}(u)$, then
$(w\wedge \sup u)\vee \inf u\in H'_{E,(\theta)}(u)$, so that $H'_{E,(\theta)}(u)$ is nonempty. \enpr

\bigskip
Let now

$$\Theta(u)=\big\{ \theta>0: \ \ \exists \, w\in H_{E,(\theta)}(u): 
\osc(\theta w\circ\psi_i)\le \osc(u)\ \forall\, i=1,...,k\big\}.$$

\bigskip
{\bf Lemma 2.} {\sl If $u,v\in\rvo$, $u$ nonconstant, and either a) or b) holds where

\quad a) $v(P)=1$ for $P$ such that $u(P)=\min u$ and $v(P)=0$ for the other $P$,

\quad b) $v(P)=-1$ for $P$ such that $u(P)=\max u$ and $v(P)=0$ for the other $P$,

then ${d\over dt}^+|_{t=0} E(u+tv)< 0$.}
\smad Fix $\bar t>0$ such that $(u+\bar t v)(P)\in [\min u,\max u]$ for every $P$. 
Suppose e.g., a holds, and let 
$A=\{P\in\vo:U(P)=\min u\}$. Then, we have
$u+\bar t v=(u\wedge\max u)\vee(\min u+\bar t)$, thus, by iv) 
we have $E(u+\bar tv)< E(u)$, As $E$ is convex, then
${d\over dt}^+|_{t=0} E(u+tv)\le {E(u+\bar t v)-E(u)\over \bar t}<0$. \enpr

\bigskip
{\bf Lemma 3.} {\sl If $u$ is nonconstant, then $H_{E,(1)}(u)\subseteq\cl''(u)$.}
\smad If $v\in\cl(u)\setminus\cl'(u)$, then 
$\tilde v:= (v\wedge\max u)\vee\min u\in \cl'(u)$ and, by iv) satisfies
$S_{(1)}(E)(\tilde v)<S_{(1)}(E)(v)$. Hence, it suffices to prove that, 
if $v\in\cl'''(u)$, then there exists $w\in\cl(u)$ such that 
$S_{(1)}(E)(w)<S_{(1)}(E)(v)$. By the definition of $\cl'''(u)$, 
there exists $\bar i$ such that $v\circ \psi_{\bar i}$ attains both 
$\max u$ and $\min u$, thus at least one of them only at points 
in $\von\setminus\vo$. Suppose for example
$v\circ \psi_{\bar i}$ attains $\min u$  only at points in $\von\setminus\vo$. Let

$$w_t:=v+t\chi_B,\quad  B:=\{Q\in\von\setminus\vo: v(Q)=\min u\},$$

for $t$ positive and sufficiently small. Then,

$$w_t\circ\psi_i=v\circ\psi_i + tr,\quad r(P_j):=\cases
1&\ \text{if}\ j\ne i, v\big(\psi_i(P_j)\big)=\min u\\ 0&\ \text{otherwise.}\endcases$$

By Lemma 2, we have ${d\over dt}^+|_{t=0} E(w_t\circ \psi_i)< 0$ if $v\circ\psi_i(P_j)=\min u$ 
holds for some $j$ including $j=i$, but not for all $j$, by v) we have
${d\over dt}^+|_{t=0} E(w_t\circ \psi_i)\le  0$. If $v\circ \psi_i(P_j)=\min u$ 
holds for all $j$, then by vi) we have
${d\over dt}^+|_{t=0} E(w_t\circ \psi_i)= 0$. If finally,  $v\circ \psi_i(P_j)=\min u$ 
never holds, then $r=0$, thus
${d\over dt}^+|_{t=0} E(w_t\circ \psi_i)=0$.  In conclusion, 
${d\over dt}^+|_{t=0} S_{(1)}(E)(w_t)< 0$, 
thus for some positive $t$ we have $S_{(1)}(E)(w_t)<S_{(1)}(E)(v)$. 
Moreover, $w_t\in\cl(u)$ as $\chi_B=0$ on $\vo$.
\enpr

\bigskip
{\bf Theorem 4.}  {\sl There exists $\delta>0$ such that, 
if $\theta\in [1,1+\delta]$, then the minimum in the definition
of $\Lb_{(\theta)}(E)(u)$ is attained only at functions in $\cl''(u)$, 
thus $\cl''(u)\supseteq [1,1+\delta]$.}
\smad By Lemma 3, there exists $\bar v\in\cl''(u)$ such that

$$S_{(1)}(E)(\bar v)\le\min \{ S_{(1)}(E)(v): v\in\cl'''(u)\} -\eta$$

for some $\eta>0$ (note that $\cl'''(u)$ is compact). thus, by continuity, 
there exists $\delta>0$ such that

$$\theta \bar v\in\cl''(u),\quad S_{(\theta)}(E)(\bar v)\le\min \Big\{ S_{(1)}(E)(v): v\in
\bigcup\limits_{t\in [\theta^{-1},1]} t\cl'''(u\Big)\}\quad \forall\, \theta\in [1,1+\delta]. \eqno (1)$$

On the other hand, for every $v\in\cl(u)$ we have 
$S_{(\theta)}(E)(v)\ge S_{(\theta)}(E)(\tilde v)$ where
$\tilde v:=(v\wedge \max u)\vee \min u\in\cl'(u)$. 
This, in combination with (1), concludes the proof. \enpr

\bigskip
As $\Theta(u)$ is clearly bounded, we can consider 
$\theta(u)=\sup \Theta(u)\in \Theta(u)$. Then, we have

\bigskip
{\bf Lemma 5.} {\sl There exists $w\big(:=w(u)\big)\in H_{E,(\theta(u))}(u)$ such that 
$\osc\big(\theta(u) w\circ\psi_i\big)\le \osc(u)$ for all $i$ and 
the inequality holds for at least one $i$.}
\smad By definition of $\Theta(u)$ and the fact that $\theta(u)\in\Theta(u)$, there exists $\bar w\in 
H_{E,(\theta(u))}(u)$ such that $\osc\big(\theta(u) w\circ\psi_i\big)\le \osc(u)$
 for all $i=1,...,k$. Let $s(w)=
\max\limits_{i=1,...,k} \osc(w\circ\psi_i)$. Then, clearly $s$ is a continuous function, and

$$s\big(\theta(u) \bar w\big)\le \osc(u). \eqno (2)$$

On the other hand, by definition of $\Theta(u)$, if $\theta_n>\theta(u)$ 
and $\theta_n\llim\limits_{n\to +\infty} \theta(u)$,
we have  $\theta_n\notin \Theta(u)$, thus every $w_n\in  H'_{E,(\theta_n)}(u)$
 satisfies $s(\theta_n w_n)>\osc(u)$.
Also, by compactness, we can and do assume $w_n\llim\limits_{n\to +\infty} \tilde w$. Clearly, 
$\tilde w\in  H'_{E,(\theta(u))}(u)$ and $s\big(\theta(u) \tilde w\big)\ge\osc(u)$, thus using (2) and
 the intermediate value Theorem on the convex set $ H_{E,(\theta(u))}(u)$, there exists 
 $w\in  H_{E,(\theta(u))}(u)$ such that $s\big(\theta(u) \tilde w\big)=\osc(u)$. \enpr

\bigskip
Now,  when $\ce:\bre^K\to\bre$, we define $S_{(\theta)}(\ce):\bre^K\to\bre$ in the same way as for 
$E:\rvo\to [0,+\infty[$, that is 

$$S_{(\theta)}(\ce)(v)=\sum\limits_{i=1}^k \ce(\theta v\circ\psi_i),
\quad S(\ce)(v)=S_{(\theta(v|_{\vo}))}(\ce)(v)\quad \forall\, v\in\bre^K.$$

Given $E\in\ca$, define $\ce_E:\bre^K\to [0,+\infty[$ as $\ce_E(v)= E(v|_{\vo}).$ 
Given $E\in\ca$ and 
$\ce:\bre^K\to\bre$, $u\in\rvo$, $v\in\bre^K$ put

$$\beta_{\theta,E}(u)=\cases {E(u)\over \Lb_{(\theta)}(E)(u)}&\quad \text{if}\ u\ \text{nonconstant},\\
1 &\quad \text{otherwise}.\endcases\quad \beta_E(v)=\beta_{\theta(v|_{\vo}),E}(v|_{\vo}).$$

$$\witi S_{E}(\ce)(v):=\beta_{E}(v)\,     S(\ce)(v), \ \ \witi \ce_{n,E}=\witi S^n_{E}(\ce_E).$$
%$$\witi S_{(\theta);n}(E)=\witi S^n_{(\theta)E}(\ce_E).$$

Note that  $\witi S_{E}(\ce):\bre^K\to[0,+\infty[$, so that the definition of 
$\witi  \ce_{n,E}$ makes sense.

\bigskip
{\bf Lemma 6.} {\sl Given $E\in\ca$ and $v\in\bre^K$, the sequence $\witi \ce_{n,E}(v)$ is increasing.}
\smad  Note that $\witi  \ce_{1,E}(v)=\witi S_{E}(\ce_E)(v)\ge E(v|_{\vo})= \ce_{0,E}(v)$ by definition as
 
$S_{(\theta)}(\ce_E)(v)\ge \Lb_{(\theta)}(E)(v|_{\vo})$. Next, we easily see that $\ce\mapsto \witi S_{E}(\ce)$ is increasing,
so that the Lemma follows by recursion.  \enpr

\bigskip
We now put on $\bre^K$ a metric $\tilde d_{\infty}$ defined by
$$\tilde d_{\infty}(v_1,v_2)=\cases \sup |v_1-v_2| & \text{if} \ v_1=v_2\ \text{on} \ \vo\\
1 & \text{otherwise} \endcases$$

in other words a sequence $v_n$ converges to $v$ in such a metric if $v_n$ tends to $v$ uniformly and
$v_n=v$ on $\vo$ eventually. Now, we put

$$\witi \ce_{\infty,E}=\Gamma(\bre^K,\tilde d_{\infty})\lim\limits_{n\to +\infty} \witi \ce_{n,E}$$

where the existence of such a $\Gamma$ limit is assured by Lemma 6.

\bigskip
{\bf Lemma 7.} {\sl The functional $\witi \ce_{\infty,E}$ is self-similar in the sense that
$$\witi \ce_{\infty,E}= \witi S_{E}\big( \witi \ce_{\infty,E}\big). \eqno(3)$$}
\smad We have

$$\witi \ce_{\infty,E}(v)=\Gamma(\bre^K,\tilde d_{\infty})\lim\limits_{n\to +\infty} \witi \ce_{n+1,E}$$
$$=\Gamma(\bre^K,\tilde d_{\infty})\lim\limits_{n\to +\infty}  \witi S_E\big( \witi \ce_{n,E}\big)(v)$$
$$=\Gamma(\bre^K,\tilde d_{\infty})\lim\limits_{n\to +\infty} \beta_E(v) \sum\limits_{i=1}^k 
\witi \ce_{n,E}(\theta(v) v\circ\psi_i)$$
$$\ge \beta_E(v) \sum\limits_{i=1}^k \Gamma(\bre^K,\tilde d_{\infty})\lim\limits_{n\to +\infty}
\witi\ce_{n,E}(\theta(v) v\circ\psi_i)$$
$$= \sum\limits_{i=1}^k\beta_E(v)\witi \ce_{\infty,E}(\theta(v) v\circ\psi_i)=  \witi S_{E}\big( \witi \ce_{\infty,E}\big)(v)$$

and the inequality $\ge$ holds in (3). On the other hand, if 
$v_{i,n}\llim\limits_{n\to +\infty} v\circ\psi_i$, we can and do assume $v_{i,n}=v\circ\psi_i$ on $\vo$.
We now take such $v_{i,n}$  such that $\witi \ce_{\infty,E}(\theta(v)v\circ\psi_i)=\lim\limits_{n\to +\infty}
\witi\ce_{n,E}(\theta(v)v_{i,n})$. Define now $v_n$ on $\vo $ by

$$v_n\circ\psi_i=v_{i,n}\ \text{on} \ K\eqno (4)$$

Formula (4) makes sense as, if $Q=\psi_i(P)=\psi_{i'}(P')$ with $P\ne P'$, then  $P,P'\in\vo$, thus
$v_{i,n}(P)=v(\psi_i(P))=v(\psi_{i'}(P'))=v_{i',n}(P')$. Moreover, we easily see that $v_n\llim\limits_{n\to +\infty}v$, and

$$\witi\ce_{n+1,E}(v_{n})=
\beta_E(v) \sum\limits_{i=1}^k \witi \ce_{n,E}(\theta(v) v_n\circ\psi_i)$$
$$=\beta_E(v) \sum\limits_{i=1}^k \witi \ce_{n,E}(\theta(v) v_{i,n})\llim\limits_{n\to+\infty}
\beta_E(v) \sum\limits_{i=1}^k \witi \ce_{\infty,E}(\theta(v)v\circ\psi_i)$$
$$= \witi S_{E}\big( \witi \ce_{\infty,E}\big)(v)$$

and, consequently,

$$\witi \ce_{\infty,E}(v)\le \lim\limits_{n\to+\infty} \witi\ce_{n+1,E}(v_{n})= \witi S_{E}\big( \witi \ce_{\infty,E}\big)(v)$$

and the Lemma is proved. \enpr